\documentclass[10pt]{amsart}

\usepackage{amsmath}
\usepackage{amssymb}
\usepackage{mathtools}
\usepackage{bm}
\usepackage{comment}
\usepackage{graphicx}
\usepackage{color}
\usepackage{hyperref}
\hypersetup{colorlinks=true, linkcolor=blue, citecolor=magenta, urlcolor=blue,
  pdftitle={On the bi-Krull property in positive semirings},
  pdfauthor={Aaditya Bilakanti, Felix Gotti, Amrit Kandasamy, Harold Polo}}
\usepackage{url}
\usepackage{algpseudocode}
\usepackage{fancyhdr}
\usepackage{tikz-cd}
\usepackage{stmaryrd}
\usepackage{calrsfs}
\usepackage[shortlabels]{enumitem}
\usepackage{array}
\usepackage{booktabs}
\usepackage{tabularx}

\definecolor{red}{rgb}{1,0,0}
\definecolor{felix}{rgb}{1,0,0}

\definecolor{blue}{rgb}{0,0,1}
\definecolor{harold}{rgb}{0,0,1}

\definecolor{green}{rgb}{0,1,0}
\definecolor{aaditya}{rgb}{0,0,1}

\definecolor{clear_blue}{rgb}{0,.5,0}
\definecolor{amrit}{rgb}{0,.5,0}

\theoremstyle{plain}
\newtheorem{theorem}{Theorem}[section]
\newtheorem{prop}[theorem]{Proposition}
\newtheorem{lemma}[theorem]{Lemma}
\newtheorem{corollary}[theorem]{Corollary}
\newtheorem*{maintheorem}{Main Theorem}

\theoremstyle{definition}
\newtheorem{example}[theorem]{Example}
\newtheorem{question}[theorem]{Question}
\theoremstyle{remark}

\numberwithin{equation}{section}

\newcommand{\nn}{\mathbb{N}}
\newcommand{\zz}{\mathbb{Z}}
\newcommand{\qq}{\mathbb{Q}}
\newcommand{\rr}{\mathbb{R}}
\newcommand{\cc}{\mathbb{C}}
\newcommand{\gp}{\mathcal{G}}
\newcommand{\F}{\mathcal{F}}
\newcommand{\A}{\mathcal{A}}
\newcommand{\atoms}{\A}

\newcommand{\red}{\mathrm{red}}
\newcommand{\Cl}{\mathcal{C}}
\DeclareMathOperator{\supp}{supp}
\DeclareMathOperator{\Tr}{Tr}
\providecommand\ldb{\llbracket}
\providecommand\rdb{\rrbracket}

\keywords{Positive semiring, Krull monoid, bi-Krull property, divisor theory,
$\omega$-primality, transfer Krull monoid}


\begin{document}
\mbox{}

\title{On the bi-Krull property in positive semirings}

\author{Aaditya Bilakanti}
\address{CMI\\MIT\\Cambridge, MA 02139}
\email{abilak@gmail.com}

\author{Felix Gotti}
\address{Department of Mathematics\\MIT\\Cambridge, MA 02139}
\email{fgotti@mit.edu}

\author{Amrit Kandasamy}
\address{CMI\\MIT\\Cambridge, MA 02139}
\email{akandas5@asu.edu}


\author{Harold Polo}
\address{Department of Mathematical Sciences, Clemson University, Clemson, SC 29634}
\email{jimenezclases@gmail.com}

\date{October 29, 2026}

\begin{abstract}
The Bi-UF Positive Conjecture, posed by Baeth, Chapman, and Gotti and recently proved by Graia, states that $\mathbb{N}_0$ is the only subsemiring of the nonnegative real numbers whose additive and multiplicative monoids are both factorial. We extend this result by characterizing the subsemirings of integral domains whose additive and multiplicative monoids are both Krull. In particular, we prove that $\mathbb{N}_0$ is the only such subsemiring of the nonnegative real numbers. We also show that the Krull condition cannot be weakened to transfer Krull.
\end{abstract}

\maketitle

\bigskip
\section{Introduction}
\label{sec:introduction}

In~\cite{BCG21}, Baeth, Chapman, and Gotti studied several classes of positive semirings (i.e., subsemirings of the nonnegative real numbers containing $1$) defined by requiring a given factorization property to hold in both the additive and the multiplicative monoids. Following their terminology, we call a positive semiring $S$ a \emph{bi-UFS} if both $(S,+)$ and $(S\setminus\{0\},\cdot)$ are factorial. As the additive monoid of $\nn_0$ is free on $\{1\}$ and its multiplicative monoid is free on the set of prime numbers, $\nn_0$ is a bi-UFS. Baeth, Chapman, and Gotti conjectured in~\cite[Conjecture~7.7]{BCG21} that $\nn_0$ is the only positive semiring that is a bi-UFS. This conjecture, known as the \emph{Bi-UF Positive Conjecture}, was recently settled by Graia~\cite[Theorem~4.3]{oG26}.

Here we study whether the conclusion of the Bi-UF Positive Conjecture continues to hold when the factoriality hypothesis on each of the two monoids is replaced by the Krull property. We say that a positive semiring $S$ is \emph{bi-Krull} provided that both its additive and multiplicative monoids are Krull.

In his proof of the Bi-UF Positive Conjecture, Graia makes use of the free structure of the additive monoid and the primality of the atoms of the multiplicative monoid (see~\cite[Section~2 and Lemma~3.2]{oG26}). Although Krull monoids need not be factorial, they admit divisor theories into free commutative monoids, and so we can still express divisibility in terms of coordinates (see~\cite[Section~2.2 and Definition~4.1]{GZ20}). On the other hand, the atoms of a Krull monoid need not be prime, but every element of a Krull monoid has finite $\omega$-primality~\cite[Theorem~4.2(2)]{GH08}. In view of these observations, we ask the following question.

\begin{question}\label{ques:bi-krull}
    If a positive semiring $S$ is bi-Krull, must $S=\nn_0$?
\end{question}

The main result of this paper is an affirmative answer to Question~\ref{ques:bi-krull}.

\begin{maintheorem}
    If a positive semiring $S$ is bi-Krull, then $S=\nn_0$.
\end{maintheorem}

We begin by showing in Proposition~\ref{prop:units} that the only multiplicative unit of a bi-Krull positive semiring is $1$. To do so, we consider the action of the multiplicative units on the coordinates afforded by a divisor theory of the additive monoid and use a polynomial of Kainrath. In Section~\ref{sec:omega}, we generalize a result of Graia~\cite[Lemma~3.2]{oG26} concerning the primality of certain elements of the multiplicative monoid. More specifically, we show that if the additive monoid of a positive semiring $S$ is Krull and a multiplicative nonunit $p\in S\setminus\{0\}$ satisfies
\[
    1+py=pz
\]
for some $y,z\in S$, then $p$ has infinite $\omega$-primality in the multiplicative monoid of $S$. Hence no such relation can hold in a bi-Krull positive semiring.

In Section~\ref{sec:algebraic}, we prove that the only algebraic numbers that can belong to a bi-Krull positive semiring are the nonnegative integers (Theorem~\ref{thm:algebraic-rigidity}). For an algebraic element $s$ of a bi-Krull positive semiring $S$, we use a divisor theory of $(S,+)$ to realize the additive monoid of $S\cap\mathbb{Q}(s)$ as the monoid of lattice points of a rational polyhedral cone. If $s$ is irrational, we use properties of this cone together with Dirichlet's unit theorem to obtain a relation of the form ruled out in Section~\ref{sec:omega}.

It remains to show that a bi-Krull positive semiring cannot contain any transcendental numbers. In Section~\ref{sec:transcendental}, we return to the polynomial of Kainrath used in the proof of Proposition~\ref{prop:units}. By homogenizing this polynomial and using the root closedness of the multiplicative monoid, we show that the support of the additive divisor vector of $s^2$ is contained in the union of the supports of the divisor vectors of $1$ and $s$. It follows that the additive divisor vectors of the elements
\[
    s,s^2,s^4,s^8,\ldots
\]
are supported on a fixed finite set. These vectors are therefore linearly dependent over $\mathbb{Q}$, from which we obtain a nonzero polynomial vanishing at $s$. This proves that every element of a bi-Krull positive semiring is algebraic (Theorem~\ref{thm:every-element-algebraic}) and so completes the proof of the Main Theorem. As a consequence of the Main Theorem, we also determine the semidomains whose additive and multiplicative monoids are both Krull: up to isomorphism, they are precisely $\nn_0$ and the Krull domains (Corollary~\ref{cor:general-semidomains}).

In the last section of this paper (Section~\ref{sec:hypotheses}), we investigate whether the conclusion of the Main Theorem still holds if one of the two monoids is only assumed to be transfer Krull. We construct a positive semiring different from $\nn_0$ whose multiplicative monoid is Krull and whose additive monoid is both root-closed and half-factorial (and, in particular, transfer Krull). However, we have not been able to decide whether $\nn_0$ is the only positive semiring whose additive monoid is Krull and whose multiplicative monoid is transfer Krull.

\bigskip
\section{Background}
\label{sec:background}

We let $\nn$, $\nn_0$, $\zz$, $\qq$, $\rr$, and $\cc$ denote the sets of positive integers, nonnegative integers, integers, rational numbers, real numbers, and complex numbers respectively. For $m, n \in \zz$ with $m \le n$, we set $\ldb m, n \rdb \coloneqq \{k \in \zz : m \le k \le n\}$. For a set $P$, $\zz^{(P)}$ denotes the group of finitely supported integer vectors indexed by $P$. We use similar notation over $\nn_0$, $\qq$, and $\rr$. The \emph{support} $\supp(v)$ of a vector $v$ is the set of coordinates at which it is nonzero.
\subsection{Krull monoids and divisor theories}

Throughout this paper, monoids are commutative and cancellative. Let $M$ be a multiplicatively written monoid. Its \emph{Grothendieck group} is denoted by $\gp(M)$, and its group of units is denoted by $M^\times$. The monoid $M$ is \emph{reduced} if $M^\times = \{1\}$. We write $M_{\red} \coloneqq M/M^\times$ for its associated reduced monoid. For $a, b \in M$, we write $a \mid_M b$ if $b = ac$ for some $c \in M$. A nonunit $a \in M$ is an \emph{atom} if $a=bc$ with $b,c\in M$ implies that $b$ or $c$ is a unit. We denote the set of atoms of $M$ by $\atoms(M)$. The monoid $M$ is \emph{atomic} if every nonunit is a finite product of atoms, and it is \emph{factorial} if these factorizations are unique up to order and associates.

Let $\F(P)$ denote the free commutative monoid on a set $P$. We identify $\F(P)$ with the additive monoid $\nn_0^{(P)}$ when working with divisor vectors. A homomorphism $\partial \colon M \to \F(P)$ is a \emph{divisor homomorphism} if $\partial(a) \le \partial(b)$ implies that $a \mid_M b$ for all $a,b \in M$, where the inequality is coordinatewise. It is a \emph{divisor theory} if, in addition, every vector in $\nn_0^{(P)}$ is the coordinatewise greatest common divisor of a finite nonempty family of vectors in $\partial(M)$.

A monoid is \emph{$v$-noetherian} if it satisfies the ascending chain condition on divisorial ideals. We use the divisor-theoretic characterization of a \emph{Krull monoid}: a monoid is Krull if it admits a divisor homomorphism into a free commutative monoid. Equivalently, it admits a divisor theory; equivalently, it is completely integrally closed and $v$-noetherian. Background on these characterizations and on divisor theories can be found in~\cite[Chapter~2]{GHK06} and~\cite[Section~4]{GZ20}. Every Krull monoid is atomic and \emph{root-closed}: if $x \in \gp(M)$ and $x^n \in M$ for some $n \in \nn$ then $x \in M$.

Recall that a monoid $M$ is \emph{completely integrally closed} if, whenever $x \in \gp(M)$ and there exists $c \in M$ such that $cx^n \in M$ for every $n \in \nn$, one has $x \in M$. 

Divisor theories are unique up to isomorphism. The saturation property of a divisor homomorphism gives
\begin{equation*}
    \partial(M) = \gp(\partial(M)) \cap \nn_0^{(P)}.
\end{equation*}
A divisor homomorphism is injective when $M$ is reduced.

For $a \in M$, the \emph{$\omega$-primality} $\omega_M(a)$ is the least $N \in \nn_0$ with the following property: whenever $a \mid_M b_1 \cdots b_m$ for $b_1, \dots, b_m \in M$, there is a subset $I \subseteq \ldb 1, m \rdb$ with $|I| \le N$ such that $a \mid_M \prod_{i \in I} b_i$. The empty product is $1$. If no such $N$ exists then $\omega_M(a) \coloneqq \infty$. We write $\omega(a)$ when the underlying monoid is clear from the context. A nonunit is prime if and only if its $\omega$-primality is $1$. Every element of a $v$-noetherian monoid has finite $\omega$-primality~\cite[Theorem~4.2]{GH08}. For a Krull monoid, a divisor theory gives the explicit bound $\omega_M(a) \le |\partial(a)|$,
where $|\partial(a)|$ is the sum of the coordinates of $\partial(a)$. Indeed, for each coordinate we can choose at most as many factors as the required multiplicity so that their total contribution reaches that multiplicity. The union of these choices has at most $|\partial(a)|$ factors, and the divisor-homomorphism property gives the divisibility.
\subsection{Positive semirings and additive divisor coordinates}

A \emph{positive semiring}, also called a \emph{positive semidomain}, is a subsemiring $S$ of $\rr_{\ge 0}$ containing $0$ and $1$. Throughout this paper, we set $S^\bullet\coloneqq S \setminus\{0\}$ and $\gp_{+}(S) \coloneqq \gp((S,+)).$
The group $\gp_{+}(S)$ identifies with the subring $S - S$ of $\rr$. We denote the set of additive atoms of $S$ by $\atoms_+(S)$. We refer to the atoms of the nonzero multiplicative monoid of a semiring as \emph{irreducibles}.

Whenever $(S,+)$ is Krull, fix a divisor theory $\partial_+ \colon (S,+) \to \nn_0^{(P)}$. Since $(S,+)$ is reduced, this map extends to an injective homomorphism $\partial_+ \colon \gp_{+}(S) \to \zz^{(P)}$. In particular,
\begin{equation}\label{eq:additive-saturation}
    S = \left\{x \in \gp_{+}(S) : \partial_+(x) \in \nn_0^{(P)}\right\}.
\end{equation}
For $\rho \in P$, let $\partial_{+, \rho}(x)$ denote the $\rho$-coordinate of $\partial_+(x)$. Outside Section~\ref{sec:hypotheses}, the set $P$ indexes the additive divisor coordinates.

\begin{lemma}\label{lem:additive-atoms}
    Let $S$ be a positive semiring whose additive monoid is Krull. Then the following statements hold.
    \begin{enumerate}
        \item[(1)] $1 \in \atoms_+(S)$.
        \smallskip

        \item[(2)] The equality $S \cap \qq = \nn_0$ holds. In particular, the element $1$ is the only rational additive atom of $S$.
        \smallskip

    \end{enumerate}
\end{lemma}

\begin{proof}
    (1) Since $(S,+)$ is nontrivial, reduced, and atomic, it has an atom $a$. If $1 = s + t$ for $s, t \in S$ then $a = as + at$. The irreducibility of $a$ implies that one of $as$ and $at$ is zero. Since $a > 0$, one of $s$ and $t$ is zero. Hence $1$ is an additive atom.
    \smallskip

    (2) Take $m/n \in S \cap \qq_{>0}$, where $m, n \in \nn$ and $\gcd(m, n) = 1$. There exist $r, t \in \zz$ such that $rm + tn = 1$, which implies $1/n = r(m/n) + t \in \gp_{+}(S)$. Since $1 \in S$, additive root closedness gives $1/n \in S$. Part~(1) now forces $n = 1$. The reverse inclusion $\nn_0 \subseteq S$ follows from the semiring axioms.
\end{proof}

\subsection{Multiplicative units}

The proof of Proposition~\ref{prop:units} uses a polynomial of Kainrath with a negative coefficient whose square has nonnegative coefficients. This polynomial appears in~\cite[Example~5.2]{CF19} as a witness to the failure of root-closedness of the nonzero multiplicative monoid of $\nn_0[x]$. Set
\[
    K(x) \coloneqq x^4 + 2x^3 - x^2 + 4x + 2.
\]
Then
\begin{equation}\label{eq:kainrath-multiple}
    (x + 1)K(x) = x^5 + 3x^4 + x^3 + 3x^2 + 6x + 2
\end{equation}
and
\begin{equation}\label{eq:kainrath-square}
    K(x)^2 = x^8 + 4x^7 + 2x^6 + 4x^5 + 21x^4 + 12x^2 + 16x + 4.
\end{equation}
In particular $K(x) > 0$ for every $x \geq 0$.

\begin{prop}\label{prop:units}
    Let $S$ be a positive semiring. If $(S,+)$ is Krull and $(S\setminus\{0\},\cdot)$ is root-closed then $S$ is multiplicatively reduced.
\end{prop}

\begin{proof}
    Fix $u \in S^{\times}$. Multiplication by $u$ is an automorphism of $(S,+)$. By uniqueness of divisor theories, it induces a permutation $\sigma$ of $P$ such that $\partial_+(ux) = \sigma(\partial_+(x))$ for every $x \in \gp_{+}(S)$, where $\sigma$ acts on divisor vectors by permuting coordinates. Set $d \coloneqq \partial_+(1)$ and $D \coloneqq \supp(d)$.

    Suppose first that every $\sigma$-orbit meeting $D$ is finite. Since $D$ is finite, there exists $N \in \nn$ such that $\sigma^N(d) = d$. Thus $\partial_+(u^N) = \partial_+(1)$, and the injectivity of $\partial_+$ gives $u^N = 1$. As $u$ is positive, we have that $u = 1$.

    Now suppose that $D$ meets an infinite $\sigma$-orbit, and choose $\rho \in D$ in that orbit. Since $D$ is finite, there exists $N \in \nn$ such that $\sigma^{jN}(\rho) \notin D$ for every $j \in \{-2,-1,1,2\}$. Thus the coordinate $\sigma^{2N}(\rho)$ belongs to the support of $\sigma^{2N}(d)$ but not to the support of $d$, $\sigma^N(d)$, $\sigma^{3N}(d)$, or $\sigma^{4N}(d)$.

    Set $w \coloneqq K(u^N) \in \gp_+(S)$. So we have
    \[
        \partial_+(w)
        = \sigma^{4N}(d) + 2\sigma^{3N}(d)
        - \sigma^{2N}(d) + 4\sigma^N(d) + 2d.
    \]
    The $\sigma^{2N}(\rho)$-coordinate of $\partial_+(w)$ is negative, so $w \notin S$. On the other hand, $w > 0$, and Equations~\eqref{eq:kainrath-multiple} and~\eqref{eq:kainrath-square} show that $w(1+u^N) \in S$ and $w^2 \in S$. The first inclusion shows that $w \in \mathcal{G}(S)$. The second inclusion and root-closedness therefore imply that $w \in S$, a contradiction.
\end{proof}

\bigskip
\section{The infinite \texorpdfstring{$\omega$}{omega}-primality}
\label{sec:omega}

In this section, we show that if $(S,+)$ is Krull and the reciprocal of a multiplicative nonunit $p$ lies in the additive group $\gp_{+}(S)$ then $\omega(p)=\infty$. As in~\cite[Lemma~3.2]{oG26}, we express the reciprocal condition as a relation $1+py=pz$. Using additive divisor coordinates, we construct products of arbitrarily many factors that are divisible by $p$, although none of their proper subproducts is. This rules out such a relation when $S$ is bi-Krull.

\begin{prop}\label{prop:omega-obstruction}
    Let $S$ be a positive semiring such that $(S,+)$ is Krull. If $p \in S^{\bullet}\setminus S^{\times}$ satisfies $1 + py = pz$ for some $y, z \in S$ then $\omega(p) = \infty$. In particular, this can not happen when $S$ is bi-Krull.
\end{prop}
\begin{proof}
    The relation $1+py=pz$ gives $1/p=z-y\in\gp_+(S)$. Since $p$ is a nonunit, $1/p\notin S$. Thus $I\coloneqq\{\rho\in P:\partial_{+,\rho}(1/p)<0\}$
    is a nonempty finite set. For every $\rho\in I$, we know that $1/p=z-y$ also gives $\partial_{+,\rho}(y)>0$. For $m\in\nn$ and $\rho\in I$, consider the polynomial
    \[
        \partial_{+,\rho}(1/p)
        +\sum_{i=1}^{m}\binom mi t^i
          \partial_{+,\rho}(p^{i-1}y^i).
    \]
    Its constant term is negative, its coefficient of $t$ is positive, and its other coefficients are nonnegative. It is therefore strictly increasing for $t\ge0$ and has exactly one positive root. That root is at most $1$. At $t=1$, the polynomial is the $\rho$-coordinate of
    \[
        \frac{(1+py)^m}{p}=p^{m-1}z^m\in S.
    \]
    For every $t>0$, increasing $m$ strictly increases the value of the polynomial, since the coefficient of $t$ strictly increases and no other coefficient decreases.

    Let $t_m$ be the largest of these roots as $\rho$ ranges over $I$. Then $0<t_{m+1}<t_m\le1.$
    Fix $m_0\ge2$, choose $K,N\in\nn$ such that $t_{m_0}<K/N<t_{m_0-1},$ and set $r\coloneqq N+Kpy\in S$. For every $j\in\nn$,
    \[
        \frac{r^j}{p}
        =N^j\frac1p+
          \sum_{i=1}^{j}\binom ji N^{j-i}K^i p^{i-1}y^i
        \in\gp_+(S).
    \]
    Its coordinates outside $I$ are nonnegative. For a coordinate in $I$, the above is $N^j$ times the corresponding polynomial at $K/N$. So all its coordinates are nonnegative if and only if when $K/N\ge t_j$. By~\eqref{eq:additive-saturation}, $p$ divides $r^{m_0}$ in $S^\bullet$, but it divides none of $r,r^2,\dots,r^{m_0-1}$. It does not divide $1$ either since it's a nonunit. So the product of $m_0$ copies of $r$ is divisible by $p$, while no proper subproduct is. Since $m_0$ was arbitrary, $\omega(p)=\infty$. This is impossible when $S^\bullet$ is Krull.
\end{proof}

\begin{corollary}\label{cor:reciprocal}
    Let $S$ be a bi-Krull positive semiring. If $p \in S$ and $1/p \in \gp_{+}(S)$ then $p = 1$.
\end{corollary}

\begin{proof}
    Suppose, by way of contradiction, that $p\ne1$. By Proposition~\ref{prop:units}, the element $p$ is a multiplicative nonunit. Write $1/p=z-y$ with $y,z\in S$. If $y=0$ then $1/p\in S$, contrary to $p$ being a nonunit. Otherwise $y,z\in S$, and $1+py=pz$ contradicts Proposition~\ref{prop:omega-obstruction}.
\end{proof}

\bigskip

\section{Algebraic elements}
\label{sec:algebraic}

By part~(2) of Lemma~\ref{lem:additive-atoms}, the equality $S\cap\qq=\nn_0$ holds for every positive semiring $S$ whose additive monoid is Krull. Here we show that if $S^\bullet$ is also Krull, then $S$ cannot contain any irrational algebraic numbers. Recall from Corollary~\ref{cor:reciprocal} that $1$ is the only unit of the ring $S-S$ belonging to a bi-Krull positive semiring $S$. If $S$ contains an irrational algebraic number $a$, we will see that the subring of $\qq(a)$ generated by $S\cap\qq(a)$ is an order in $\qq(a)$ and so has infinitely many units by Dirichlet's unit theorem. We use a divisor theory of $(S,+)$ to describe which elements of this order belong to $S$ and, in particular, to find a unit other than $1$ in $S$. Let $\overline{\qq}$ denote the field of algebraic numbers in $\cc$.

\begin{theorem}\label{thm:algebraic-rigidity}
    If $S$ is a bi-Krull positive semiring, then $S\cap\overline{\qq}=\nn_0$.
\end{theorem}
\nopagebreak[4]
\begin{proof}
    Fix an algebraic element $a\in S$, and set $K\coloneqq\qq(a)$ and $d\coloneqq[K:\qq]$. Suppose, towards a contradiction, that $d>1$. We let $S_K\coloneqq S\cap K$ and $L\coloneqq S_K-S_K$. If $b,c\in S_K$ and $\partial_+(b)\leq\partial_+(c)$ then $c-b\in S_K$ by~\eqref{eq:additive-saturation}. Thus the restriction of $\partial_+$ to $S_K$ is a divisor homomorphism, and so $(S_K,+)$ is Krull. Because $L$ embeds in the free abelian group $\zz^{(P)}$, its additive group is free abelian. Furthermore, because $L\subseteq K$ and $1,a,\dots,a^{d-1}\in L$, its rank is $d$.

    For each $b\in L$, multiplication by $b$ is represented by a $d\times d$ integer matrix with respect to a $\zz$-basis of $L$. Evaluating the identity of multiplication maps obtained from the Cayley--Hamilton theorem at $1$, we see that $b$ satisfies a monic polynomial in $\zz[x]$. Hence $L$ is a subring of finite index in the ring of integers of $K$, that is, an order in $K$.

    We now describe $S_K$ in terms of the divisor coordinates of $L$. Let $P_0$ be the union of the supports of the divisor vectors of a fixed $\zz$-basis of $L$. As $P_0$ is finite and $\partial_+(L)\subseteq\zz^{P_0}$, we can extend the restriction of $\partial_+$ to $L$ to a linear map from $V\coloneqq L\otimes_{\zz}\rr$ to $\rr^{P_0}$, which we also denote by $\partial_+$. This extension is injective because its matrix has integer entries and rank $d$. Consider the cone
    \[
        C\coloneqq\{x\in V:\partial_{+,\rho}(x)\geq 0\text{ for every }\rho\in P_0\}.
    \]
    Since $C$ is defined by finitely many homogeneous linear inequalities with integer coefficients, it is a closed rational polyhedral cone. Also, $C\cap(-C)=\{0\}$ because $\partial_+$ is injective, and $C$ has nonempty interior because it contains the basis $1,a,\dots,a^{d-1}$ of $V$. By~\eqref{eq:additive-saturation}, we have $C\cap L=S_K$.

    The rational polyhedral cone $C$ is generated by finitely many rational vectors (see~\cite[Theorem~6]{Kaibel11}). After multiplying these vectors by suitable positive integers, we may assume that they belong to $L$, and therefore to $S_K$. Thus every element of $C$ is a nonnegative real linear combination of elements of $S_K$. The product of two such combinations is another such combination because $S_K$ is a semiring. Hence $C$ is closed under multiplication.

    Because $L$ spans $K$ over $\qq$, we can identify $V$ with $K\otimes_{\qq}\rr$. Let $\sigma_0,\dots,\sigma_{r-1}$ consist of the real embeddings of $K$ and one embedding from each pair of nonreal conjugate embeddings, where $\sigma_0$ is the inclusion of $K$ into $\rr$. These embeddings allow us to identify $V$ with a product of copies of $\rr$ and $\cc$, in which addition and multiplication are performed coordinatewise. We use $\sigma_i$ also for the corresponding coordinate maps on $V$, with the remaining embeddings of $K$ extending to the complex conjugates of these maps. Set $e_0\coloneqq(1,0,\dots,0)\in V$. As we will see, there is a unit $u$ of $L$ with $\sigma_0(u)>1$ such that $u^n/\sigma_0(u)^n\to e_0$. To show that some power of $u$ belongs to $C\cap L$, it suffices to show that $e_0$ belongs to the interior $C^\circ$ of $C$. We establish this in three claims.
    \smallskip

    \noindent\refstepcounter{claim}\label{claim:positive-embedding}%
    \textsc{Claim~\theclaim.} For every $x\in C\setminus\{0\}$, the inequality $\sigma_0(x)>0$ holds.
    \par\nobreak\smallskip
    \noindent\textsc{Proof of Claim~\ref{claim:positive-embedding}.}
    Every nonzero element of $C$ is a nonnegative linear combination of nonzero elements of $S_K$ with at least one positive coefficient. Since $\sigma_0$ is linear and positive on $S_K\setminus\{0\}$, it is also positive on $C\setminus\{0\}$.
    \smallskip

    \noindent\refstepcounter{claim}\label{claim:dominant-embedding}%
    \textsc{Claim~\theclaim.} For every embedding $\sigma\colon K\to\cc$ and every $x\in C$,
    \begin{equation}\label{eq:dominant-embedding}
        |\sigma(x)|\leq\sigma_0(x).
    \end{equation}
    \noindent\textsc{Proof of Claim~\ref{claim:dominant-embedding}.}
    Suppose that~\eqref{eq:dominant-embedding} fails for some embedding of $K$ and some $x\in C$. Set $R\coloneqq\max\{|\sigma_i(x)|:i\in\ldb 0,r-1\rdb\}$, and note that $R>\sigma_0(x)\geq 0$. Since $C$ is closed under multiplication, every term of the sequence $(x^n/R^n)_{n\geq 1}$ belongs to $C$. Each coordinate of $x^n/R^n$ has absolute value at most $1$, so the sequence is bounded and has a convergent subsequence whose limit $z$ belongs to $C$, as $C$ is closed. Moreover, $|\sigma_i(x)|=R$ for some $i$, and so $|\sigma_i(z)|=1$ and $z\neq 0$. On the other hand, $\sigma_0(x)<R$ implies that $\sigma_0(z)=0$, contrary to Claim~\ref{claim:positive-embedding}.
    \smallskip

    \noindent\refstepcounter{claim}\label{claim:interior-idempotent}%
    \textsc{Claim~\theclaim.} The element $e_0$ belongs to $C^\circ$.
    \par\nobreak\smallskip
    \noindent\textsc{Proof of Claim~\ref{claim:interior-idempotent}.}
    For each $i\in\ldb 1,r-1\rdb$, the real polynomial $|\sigma_i(x)|^2-\sigma_0(x)^2$ is nonzero because it takes the value $-1$ at $e_0$. The product of these polynomials cannot vanish on a nonempty open set, so there is an $x\in C^\circ$ at which none of them vanishes. It follows from Claims~\ref{claim:positive-embedding} and~\ref{claim:dominant-embedding} that $\sigma_0(x)>0$ and $|\sigma_i(x)|<\sigma_0(x)$ for every $i\in\ldb 1,r-1\rdb$. Therefore
    \[
        \lim_{n\to\infty}\frac{x^n}{\sigma_0(x)^n}=e_0.
    \]
    Since $C$ is closed and contains every term of this sequence, $e_0\in C$.

    We still need to show that $\partial_{+,\rho}(e_0)>0$ for each $\rho\in P_0$. Recall that the trace $\Tr_{K/\qq}(b)$ of $b\in K$ is the sum of its images under all $d$ embeddings of $K$ into $\cc$. If $c\in K\setminus\{0\}$ then the $\qq$-linear map $y\mapsto\Tr_{K/\qq}(cy)$ takes the nonzero value $d$ at $c^{-1}$. The assignment of this map to $c$ is therefore an injective linear map from $K$ to the space of $\qq$-linear maps from $K$ to $\qq$. Both spaces have dimension $d$, so every such map has this form.

    Now fix $\rho\in P_0$. By the choice of $P_0$, the restriction of $\partial_{+,\rho}$ to $K$ is a nonzero $\qq$-linear map to $\qq$. Hence $\partial_{+,\rho}(y)=\Tr_{K/\qq}(cy)$ for every $y\in K$ and some $c\in K\setminus\{0\}$. This equality extends linearly to $V$, where the trace is the sum of the real coordinates and twice the sum of the real parts of the complex coordinates. In particular,
    \[
        \partial_{+,\rho}(e_0)=\Tr_{K/\qq}(ce_0)=\sigma_0(c)\neq 0.
    \]
    Since $e_0\in C$, we have $\partial_{+,\rho}(e_0)>0$. As this holds for every $\rho$ in the finite set $P_0$, we conclude that $e_0\in C^\circ$.
    \medskip

    It remains to find a unit of $L$ other than $1$ in $S_K$. Because $d>1$ and $K$ has a real embedding, $r\geq 2$. By Dirichlet's unit theorem for orders (see~\cite[Section~3]{kC}), the vectors $(\log|\sigma_i(u)|)_{i=0}^{r-1}$, as $u$ ranges over $L^\times$, form a lattice spanning the hyperplane
    \[
        \bigg\{(t_0,\dots,t_{r-1})\in\rr^r:\sum_{i=0}^{r-1}w_it_i=0\bigg\},
    \]
    where $w_i=1$ if $\sigma_i$ is real and $w_i=2$ otherwise. Since $w_0=1$ and $\sum_{i=0}^{r-1}w_i=d$, the vector $(d-1,-1,\dots,-1)$ belongs to this hyperplane. Each point of the hyperplane is within a fixed distance of a point of the lattice. Therefore, by considering a sufficiently large positive multiple of $(d-1,-1,\dots,-1)$, we can find $u\in L^\times$ such that $\log|\sigma_0(u)|>\log|\sigma_i(u)|$ for every $i\in\ldb 1,r-1\rdb$. The weighted sum of these logarithms is zero, so $\log|\sigma_0(u)|>0$.

    Replacing $u$ by its square, we obtain $\sigma_0(u)>1$ and $|\sigma_i(u)|<\sigma_0(u)$ for every $i\in\ldb 1,r-1\rdb$. In particular, the sequence $(u^m/\sigma_0(u)^m)_{m\geq 1}$ converges to $e_0$. By Claim~\ref{claim:interior-idempotent}, we have $e_0\in C^\circ$, and so $u^m/\sigma_0(u)^m\in C$ for all sufficiently large $m\in\nn$. Since $C$ is a cone, it follows that $u^m\in C\cap L=S_K$ for every sufficiently large $m\in\nn$. However, $u^m>1$ and $u^{-m}\in L\subseteq\gp_+(S)$, contrary to Corollary~\ref{cor:reciprocal}. Thus $d=1$, and it follows from Lemma~\ref{lem:additive-atoms}(2) that $a\in\nn_0$.
\end{proof}

\bigskip

\section{Bi-Krull positive semirings}
\label{sec:transcendental}

By Theorem~\ref{thm:algebraic-rigidity}, it remains to show that every element of a bi-Krull positive semiring is algebraic. If $s$ is a transcendental element of such a semiring, then $s,s^2,s^4,\dots$ are linearly independent over $\qq$. Since $\partial_+$ is additive and injective, their divisor vectors are also linearly independent. We will show, however, that all these vectors have nonzero entries in only finitely many coordinates. For this, it suffices to check that any coordinate that is zero in both $\partial_+(1)$ and $\partial_+(s)$ is also zero in $\partial_+(s^2)$.

To compare these coordinates, we return to the polynomial of Kainrath used in the proof of Proposition~\ref{prop:units} and consider the polynomial in two variables
\begin{equation}\label{eq:homogeneous-kainrath}
    K(x,y)\coloneqq x^4+2x^3y-x^2y^2+4xy^3+2y^4.
\end{equation}
The term $-x^2y^2$ will allow us to bound the coordinates of $\partial_+(s^2)$ by taking $x=s$ and letting $y$ range over $\nn$. First, we establish the following lemma.

\begin{lemma}\label{lem:homogeneous-closure}
    Let $S$ be a positive semiring such that $S^\bullet$ is root-closed. Then $K(x,y)\in S$ for all $x,y\in S$.
\end{lemma}

\begin{proof}
    The case $xy=0$ is immediate from~\eqref{eq:homogeneous-kainrath}. Suppose that $x,y\in S^\bullet$. Then $K(x,y)=y^4K(x/y)>0$ because $K(x)$ is positive on $\rr_{\geq 0}$. By~\eqref{eq:kainrath-multiple} and~\eqref{eq:kainrath-square}, the polynomials $(x+1)K(x)$ and $K(x)^2$ have nonnegative integer coefficients, and the same is true of their homogenizations $(x+y)K(x,y)$ and $K(x,y)^2$. Both $(x+y)K(x,y)$ and $K(x,y)^2$ therefore belong to $S^\bullet$. Since $x+y\in S^\bullet$, we have $K(x,y)\in\gp(S^\bullet)$. The root-closedness of $S^\bullet$ now guarantees that $K(x,y)\in S$.
\end{proof}

\begin{lemma}\label{lem:support-containment}
    Let $S$ be a positive semiring such that $(S,+)$ is Krull and $S^\bullet$ is root-closed. Then, for each $s\in S$,
    \[
        \supp(\partial_+(s^2))\subseteq\supp(\partial_+(s))\cup\supp(\partial_+(1)).
    \]
\end{lemma}

\begin{proof}
    Fix $s\in S$ and $\rho\in P$ with $\partial_{+,\rho}(s)=\partial_{+,\rho}(1)=0$. We need to show that $\partial_{+,\rho}(s^2)=0$. For every $m\in\nn$, we have $K(s,m)\in S$ by Lemma~\ref{lem:homogeneous-closure}. The additivity of $\partial_+$, together with the identity $K(s,m)+m^2s^2=s^4+2ms^3+4m^3s+2m^4$, yields
    \begin{equation}\label{eq:scaled-kainrath}
        \begin{split}
            \partial_+(K(s,m))+m^2\partial_+(s^2)
                ={}&\partial_+(s^4)+2m\partial_+(s^3)\\
                   &+4m^3\partial_+(s)+2m^4\partial_+(1).
        \end{split}
    \end{equation}
    The last two terms on the right of~\eqref{eq:scaled-kainrath} have zero $\rho$-coordinate, whereas $\partial_{+,\rho}(K(s,m))\geq 0$. Consequently,
    \[
        m^2\partial_{+,\rho}(s^2)\leq\partial_{+,\rho}(s^4)+2m\partial_{+,\rho}(s^3).
    \]
    Dividing this inequality by $m^2$ and letting $m\to\infty$, we find that $\partial_{+,\rho}(s^2)=0$.
\end{proof}

The preceding lemma suffices to establish algebraicity even when $S^\bullet$ is only assumed to be root-closed.

\begin{theorem}\label{thm:every-element-algebraic}
    Let $S$ be a positive semiring such that $(S,+)$ is Krull and $S^\bullet$ is root-closed. Then every element of $S$ is algebraic.
\end{theorem}

\begin{proof}
    For a fixed $s\in S$, set $P_s\coloneqq\supp(\partial_+(s))\cup\supp(\partial_+(1))$. By induction on $k$, Lemma~\ref{lem:support-containment} implies that $\supp(\partial_+(s^{2^k}))\subseteq P_s$ for every $k\in\nn_0$. In particular, the $|P_s|+1$ vectors $\partial_+(s^{2^k})$ with $k\in\ldb 0,|P_s|\rdb$ lie in the vector space $\qq^{P_s}$ of dimension $|P_s|$ and are therefore linearly dependent. Clearing denominators, we can choose integers $c_0,\dots,c_{|P_s|}$, not all zero, such that
    \[
        \partial_+\left(\sum_{k=0}^{|P_s|}c_k s^{2^k}\right)=\sum_{k=0}^{|P_s|}c_k\partial_+(s^{2^k})=0.
    \]
    Since $\partial_+$ is injective on $\gp_+(S)$, we must have $\sum_{k=0}^{|P_s|}c_ks^{2^k}=0$. In particular, $s$ is a root of the nonzero polynomial $\sum_{k=0}^{|P_s|}c_kx^{2^k}\in\zz[x]$.
\end{proof}

We are now in a position to prove the Main Theorem stated in Section~\ref{sec:introduction}.

\begin{proof}[Proof of the Main Theorem]
    Let $S$ be a bi-Krull positive semiring. The monoid $(S,+)$ is Krull, and $S^\bullet$, being Krull, is root-closed. Hence every element of $S$ is algebraic by Theorem~\ref{thm:every-element-algebraic}. On the other hand, Theorem~\ref{thm:algebraic-rigidity} ensures that $S$ contains no algebraic numbers outside $\nn_0$. Since $\nn_0\subseteq S$, we conclude that $S=\nn_0$.
\end{proof}

The Krull hypothesis on $S^\bullet$ in the Main Theorem can be weakened. Indeed, if $(S,+)$ is Krull and $S^\bullet$ is root-closed, then Theorem~\ref{thm:every-element-algebraic} still applies. If, in addition, every element of $S^\bullet$ has finite $\omega$-primality, Corollary~\ref{cor:reciprocal} follows from Propositions~\ref{prop:units} and~\ref{prop:omega-obstruction} as before. The proof of Theorem~\ref{thm:algebraic-rigidity} uses the multiplicative Krull hypothesis only through this corollary, so we still have $S=\nn_0$. In particular, it suffices to assume that $(S,+)$ is Krull and $S^\bullet$ is root-closed and $v$-noetherian because every element of a $v$-noetherian monoid has finite $\omega$-primality~\cite[Theorem~4.2]{GH08}.

We use the Main Theorem to characterize the subsemirings of integral domains whose additive and multiplicative monoids are Krull. Recall that such a semiring $S$ is additively reduced if $S\cap(-S)=\{0\}$.

\begin{corollary}[General semidomains]\label{cor:general-semidomains}
    Let $D$ be an integral domain and let $S\subseteq D$ be a unital subsemiring. The monoids $(S,+)$ and $S^\bullet$ are Krull if and only if one of the following conditions holds:
    \begin{enumerate}
        \item[(1)] $S$ is additively reduced, $\operatorname{char}(D)=0$, and $S=\nn_0 1_D$;
        \item[(2)] $S$ is a Krull domain.
    \end{enumerate}
\end{corollary}

\begin{proof}
    Suppose first that both $(S,+)$ and $S^\bullet$ are Krull. If $(S,+)$ is reduced then $D$ has characteristic zero (otherwise $-1_D$ would belong to $S$). In this case, the proof of Lemma~\ref{lem:homogeneous-closure} still applies: for nonzero $x,y\in S$, neither $(x+y)K(x,y)$ nor $K(x,y)^2$ can vanish because both are sums of nonzero elements of $S$. The proofs of Lemma~\ref{lem:support-containment} and Theorem~\ref{thm:every-element-algebraic} then show that every element of $S$ is algebraic over $\qq$.

    Each element of the ring $R\coloneqq S-S$ is therefore algebraic over $\qq$. The relation $a\leq b$ if $b-a\in S$ defines a partial order on $R$. Since $R=S-S$, any two elements of $R$ have a common upper bound, and so the order is directed. By Zorn's lemma, this order extends to a maximal partial order on $R$, which is still directed. By~\cite[Theorem~1]{Ma04}, the latter order is Archimedean and total. Hence $R$ admits an order-preserving embedding into $\rr$. The image of $S$ under this embedding is a bi-Krull positive semiring, and we conclude from the Main Theorem that $S=\nn_0 1_D$.

    If $(S,+)$ is not reduced, there exists $u\in S^\bullet$ such that $-u\in S$. Thus $-1_D=(-u)/u\in\gp(S^\bullet)$. The equality $(-1_D)^2=1_D$ and the root-closedness of $S^\bullet$ imply $-1_D\in S$, and so $S$ is a ring. Since $S^\bullet$ is Krull, $S$ is a Krull domain by~\cite[Section~3]{Krause89}.

    Conversely, both monoids of $\nn_0$ are factorial, whereas the additive monoid of a Krull domain is a group and its multiplicative monoid is Krull.
\end{proof}

\bigskip
\section{The scope of the hypotheses}
\label{sec:hypotheses}

The Main Theorem does not extend to arbitrary semidomains. Indeed, every abelian group is Krull, since the trivial homomorphism to the free monoid on the empty set is a divisor theory. Thus both monoids of a Dedekind domain contained in $\cc$ (or of a subfield of $\cc$) are Krull. See also Corollary~\ref{cor:general-semidomains} for the more general characterization.

Here we return to positive semirings and consider whether either of the two Krull hypotheses in the Main Theorem can be replaced by the transfer Krull property. We will see that the answer is negative for the additive monoid, whereas the corresponding question for the multiplicative monoid remains open.

We recall the following terminology (see~\cite[Sections~2, 4, and~5]{GZ20}). Let $M$ and $M_0$ be multiplicatively written monoids. A homomorphism $\theta\colon M\to M_0$ is a \emph{transfer homomorphism} if the following conditions hold.
\begin{enumerate}
    \item[(1)] $M_0=\theta(M)M_0^\times$ and $\theta^{-1}(M_0^\times)=M^\times$.
    \smallskip

    \item[(2)] If $u\in M$ and $\theta(u)=bc$ for some $b,c\in M_0$ then there exist $v,w\in M$ such that $u=vw$, $\theta(v)\in bM_0^\times$, and $\theta(w)\in cM_0^\times$.
    \smallskip
\end{enumerate}
When $M$ is written additively, the equality $u=vw$ in condition~(2) is replaced by $u=v+w$. A monoid is \emph{transfer Krull} if it admits a transfer homomorphism into a Krull monoid. If $\theta$ is a transfer homomorphism then an element $a\in M$ is an atom if and only if $\theta(a)$ is an atom of $M_0$. Moreover, condition~(2) allows us to lift any factorization of $\theta(a)$ into atoms to a factorization of $a$ into atoms. In particular, if $M_0$ is atomic then so is $M$, and therefore every transfer Krull monoid is atomic.

An atomic monoid is \emph{half-factorial} if any two factorizations of the same element into atoms have the same length. Clearly, every Krull monoid is transfer Krull. Every half-factorial monoid $M$ is also transfer Krull: if $M$ is not a group, the map taking each element of $M$ to the common length of its factorizations is a transfer homomorphism to $(\nn_0,+)$. Indeed, this map is surjective because $M$ contains an atom, and one can verify condition~(2) by splitting a factorization into two parts of the prescribed lengths. If $M$ is a group then it is already Krull.

\subsection{The additive hypothesis}

By part~(2) of Lemma~\ref{lem:additive-atoms}, a positive semiring whose additive monoid is Krull cannot contain any rational numbers outside $\nn_0$. We next show that the same conclusion holds under the weaker assumption that the additive monoid is transfer Krull.

\begin{prop}\label{prop:transfer-rational}
    Let $S$ be a positive semiring such that $(S,+)$ is transfer Krull. Then $1\in\atoms_+(S)$ and $S\cap\qq=\nn_0$.
\end{prop}

\begin{proof}
    Let $\theta\colon(S,+)\to M_0$ be a transfer homomorphism into a Krull monoid $M_0$. After composing $\theta$ with the canonical transfer homomorphism from $M_0$ to its associated reduced monoid, we may assume that $M_0$ is reduced. Fix a divisor theory $\partial\colon M_0\to\nn_0^{(P)}$. Since $M_0$ is atomic, $(S,+)$ is also atomic. The proof of part~(1) of Lemma~\ref{lem:additive-atoms} requires only that $(S,+)$ be reduced and atomic, so $1\in\atoms_+(S)$. In particular, $\theta(1)$ is an atom of $M_0$.

    It is clear that $\nn_0\subseteq S\cap\qq$. To prove the reverse inclusion, suppose that $m/n\in S$ for some relatively prime positive integers $m$ and $n$ with $n\ge 2$. Applying $\partial\circ\theta$ to the equality $m\cdot1=n(m/n)$ in $(S,+)$, we obtain $m\,\partial(\theta(1))=n\,\partial(\theta(m/n))$. Since $\gcd(m,n)=1$, there is a vector $D\in\nn_0^{(P)}$ such that $\partial(\theta(1))=nD$ and $\partial(\theta(m/n))=mD$. By B\'ezout's identity, $D$ is an integer linear combination of $\partial(\theta(1))$ and $\partial(\theta(m/n))$, and so belongs to $\gp(\partial(M_0))$. Therefore $D=\partial(y)$ for some $y\in M_0$ by the saturation property of $\partial$. The injectivity of $\partial$ now implies that $\theta(1)=y^n$. However, $\theta(1)$ is an atom, which is impossible because $y$ is a nonunit and $n\ge 2$. Hence $S\cap\qq\subseteq\nn_0$.
\end{proof}

The conclusion of the Main Theorem, however, does not hold if the additive monoid is only assumed to be transfer Krull. To see this, we return to the semiring introduced in~\cite[Theorem~5.5]{GGHL26}, where both of its monoids were shown to be half-factorial. In the next proposition, we determine which of the two monoids is Krull and compute the divisor class group of the multiplicative monoid. Recall that if $M$ is a Krull monoid with divisor theory $\partial\colon M\to\nn_0^{(P)}$ then its \emph{divisor class group} is $\Cl(M)\coloneqq\zz^{(P)}/\gp(\partial(M))$. Recall also that $M$ is factorial if and only if $\Cl(M)$ is trivial.

\begin{samepage}
\begin{prop}\label{prop:laurent-semiring}
    Let $t>1$ be a transcendental real number, and set
    \[
        S_t\coloneqq\{0\}\cup\{f(t):f\in\zz[x,x^{-1}],\ f(t)>0,\ f(1)\ge 1\}.
    \]
    Then $S_t$ is a positive semiring with $t\in S_t\setminus\nn_0$, and the following statements hold.
    \begin{enumerate}
        \item[(1)] The additive monoid $(S_t,+)$ is half-factorial and root-closed but is not Krull.
        \smallskip

        \item[(2)] The multiplicative monoid $S_t^\bullet$ is a half-factorial Krull monoid with group of units $\{t^k:k\in\zz\}$ and divisor class group of order two.
        \smallskip
    \end{enumerate}
    In particular, the hypothesis that $(S,+)$ is Krull in the Main Theorem, Proposition~\ref{prop:units}, and Theorem~\ref{thm:every-element-algebraic} can be replaced neither by the hypothesis that $(S,+)$ is transfer Krull nor by the hypothesis that $(S,+)$ is root-closed and half-factorial.
\end{prop}
\end{samepage}

\begin{proof}
    Set $D\coloneqq\zz[x,x^{-1}]$. Evaluation at $t$ defines an injective map from $D$ to $\rr$ because $t$ is transcendental. The sum and the product of two Laurent polynomials in $D$ taking positive values at both $1$ and $t$ also take positive values at both points. Thus $S_t$ is closed under addition and multiplication. Since $0,1,t\in S_t$ and $t\notin\nn_0$, the set $S_t$ is a positive semiring distinct from $\nn_0$.
    \smallskip

    (1) Every element $f(t)\in S_t^\bullet$ with $f(1)=1$ is an additive atom: the sum of two Laurent polynomials taking positive integer values at $1$ cannot take the value $1$ at $1$. On the other hand, if $f(1)\ge 2$ then we can choose $k\in\zz$ such that $t^k<f(t)$ because $t>1$. Then $(f-x^k)(t)>0$ and $(f-x^k)(1)=f(1)-1\ge 1$, so we can write $f(t)$ as the sum of two elements $t^k$ and $(f-x^k)(t)$ of $S_t^\bullet$. By induction on $f(1)$, each $f(t)\in S_t^\bullet$ factors additively into atoms. Moreover, each such factorization has length $f(1)$ because the Laurent polynomial corresponding to each additive atom takes the value $1$ at $1$. Therefore $(S_t,+)$ is half-factorial.

    To see that $\gp_+(S_t)=\{f(t):f\in D\}$, note that for every $f\in D$ and every sufficiently large $N\in\nn$, we can write $f(t)=(f+N)(t)-N$ as the difference of two elements of $S_t$. Now suppose that $nf(t)\in S_t^\bullet$ for some $f\in D$ and $n\in\nn$. Then $f(t)>0$ and, since $nf(1)$ is positive, the integer $f(1)$ must be positive. Hence $f(t)\in S_t$, and so $(S_t,+)$ is root-closed. However, $t-1\in\gp_+(S_t)\setminus S_t$ while $1+n(t-1)\in S_t$ for every $n\in\nn$. Thus $(S_t,+)$ is not completely integrally closed and therefore is not Krull.
    \smallskip

    (2) Since $D$ is obtained from $\zz[x]$ by inverting $x$, it is a factorial domain whose units are the monomials $\pm x^k$ with $k\in\zz$. By the injectivity of evaluation at $t$, every unit of $S_t^\bullet$ comes from a unit of $D$ taking a positive value at $1$. Conversely, $t^k$ and $t^{-k}$ belong to $S_t$ for every $k\in\zz$. Consequently, $(S_t^\bullet)^\times=\{t^k:k\in\zz\}$.

    Let $P$ be a set of representatives of the associate classes of prime elements of $D$ not vanishing at $1$, chosen so that $p(1)>0$ for every $p\in P$. Each $f\in D$ with $f(1)>0$ can be written uniquely in the form
    \[
        f=x^k\prod_{p\in P}p^{v_p},
    \]
    where $k\in\zz$ and $v\in\nn_0^{(P)}$. Let $\partial\colon S_t^\bullet\to\nn_0^{(P)}$ be the homomorphism taking $f(t)$ to its vector of prime exponents $v$. Set $P^-\coloneqq\{p\in P:p(t)<0\}$. The sign of $f(t)$ is determined by the number of prime factors in $P^-$, counted with multiplicity. Furthermore, every vector of nonnegative prime exponents occurs for some $f\in D$ with $f(1)>0$. Thus
    \begin{equation}\label{eq:laurent-divisor-image}
        \partial(S_t^\bullet)=\left\{v\in\nn_0^{(P)}:\sum_{p\in P^-}v_p\in 2\zz\right\}.
    \end{equation}
    \smallskip

    \noindent\textsc{Claim.} The homomorphism $\partial$ is a divisor theory.
    \par\nobreak\smallskip

    \noindent\textsc{Proof of Claim.}
    Suppose that $f(t),g(t)\in S_t^\bullet$ satisfy $\partial(f(t))\leq\partial(g(t))$. Then $f$ divides $g$ in $D$, so $g=fh$ for some $h\in D$. Because $g$ and $f$ take positive values at both $1$ and $t$, so does $h$. Hence $h(t)\in S_t^\bullet$, and therefore $f(t)$ divides $g(t)$ in $S_t^\bullet$. This proves that $\partial$ is a divisor homomorphism.

    To check the remaining condition in the definition of a divisor theory, observe that $P^-$ is infinite: if $a,b\in\nn$ are relatively prime with $1<a/b<t$ then $a-bx$ is a prime element of $D$ that is positive at $1$ and negative at $t$. Distinct fractions $a/b$ produce nonassociate primes. Take $v\in\nn_0^{(P)}$ not belonging to $\partial(S_t^\bullet)$, and choose distinct primes $p,q\in P^-$. By~\eqref{eq:laurent-divisor-image}, increasing either the $p$-coordinate or the $q$-coordinate of $v$ by $1$ produces a vector in $\partial(S_t^\bullet)$. The coordinatewise greatest common divisor of the two resulting vectors is $v$. Since every other vector is already in the image of $\partial$, the claim follows.
    \medskip

    Having shown that $S_t^\bullet$ is Krull, we proceed to compute its divisor class group. By~\eqref{eq:laurent-divisor-image}, every vector in $\gp(\partial(S_t^\bullet))$ has even coordinate sum over $P^-$. Conversely, any integer vector with this property is the difference of two nonnegative vectors with the same property: one need only add sufficiently large even integers to its negative coordinates. Because $P^-$ is nonempty, this subgroup has index two in $\zz^{(P)}$. Hence $\Cl(S_t^\bullet)\cong\zz/2\zz$.

    It remains to check half-factoriality. By~\eqref{eq:laurent-divisor-image} and the divisor-homomorphism property of $\partial$, the atoms of $S_t^\bullet$ are, up to associates, the elements $p(t)$ with $p\in P\setminus P^-$ and the elements $p(t)q(t)$ with $p,q\in P^-$ (allowing $p=q$). For an element $f(t)\in S_t^\bullet$ with divisor vector $v$, the length of any factorization into atoms is
    \[
        \sum_{p\in P\setminus P^-}v_p+\frac12\sum_{p\in P^-}v_p.
    \]
    This proves half-factoriality. The last assertion of the proposition follows from the fact that the transcendental number $t$ is a nontrivial unit of $S_t^\bullet$.
\end{proof}

The semiring $S_t$ in Proposition~\ref{prop:laurent-semiring} contains no algebraic numbers outside $\nn_0$. Indeed, suppose that $f\in\zz[x,x^{-1}]$ is nonconstant and $\alpha=f(t)$ is algebraic. After multiplying $f(x)-\alpha$ by a suitable power of $x$, we obtain a nonzero polynomial in $\qq(\alpha)[x]$ vanishing at $t$. Consequently, $t$ would be algebraic over $\qq(\alpha)$ and hence over $\qq$, a contradiction. Thus $S_t\cap\overline{\qq}=\nn_0$.

We do not know whether Theorem~\ref{thm:algebraic-rigidity} remains true if the additive monoid is only assumed to be transfer Krull and the multiplicative monoid is Krull.

\subsection{The multiplicative hypothesis}

In Section~\ref{sec:transcendental}, we saw that the conclusion of the Main Theorem still holds if $(S,+)$ is Krull, $S^\bullet$ is root-closed, and every element of $S^\bullet$ has finite $\omega$-primality. In particular, the multiplicative monoid need only be root-closed and $v$-noetherian. The following example shows that root closedness of the multiplicative monoid is not sufficient, even if the additive monoid is factorial.

\begin{example}\label{ex:sqrt2}
    Set $S=\nn_0[\sqrt{2}]$, and note that $(S,+)$ is free on $\{1,\sqrt{2}\}$. We first show that $S^\bullet$ is root-closed. Suppose that $x\in\gp(S^\bullet)$ satisfies $x^n=c+d\sqrt{2}$ for some $n\in\nn$ and $c,d\in\nn_0$. Then $x$ is a root of the monic polynomial $(X^n-c)^2-2d^2\in\zz[X]$. Since $x\in\qq(\sqrt{2})$ and the ring of integers of $\qq(\sqrt{2})$ is $\zz[\sqrt{2}]$, we can write $x=a+b\sqrt{2}$ for some $a,b\in\zz$. Let $\bar x$ denote the conjugate of $x$. Notice that $|\bar x|^n=|c-d\sqrt{2}|\le x^n$. Because $x>0$, we have $|\bar x|\le x$, and so both $a=(x+\bar x)/2$ and $b=(x-\bar x)/(2\sqrt{2})$ are nonnegative. Therefore $x\in S^\bullet$.

    Although $S^\bullet$ is root-closed, it contains an element of infinite $\omega$-primality. Indeed, the equality $1+(1+\sqrt{2})=(1+\sqrt{2})\sqrt{2}$ and Proposition~\ref{prop:omega-obstruction} imply that $\omega(1+\sqrt{2})=\infty$, because every nonzero element of $S$ is at least $1$. Thus $S^\bullet$ is not Krull.

    We also show that $S^\bullet$ is not half-factorial. If $1+r\sqrt{2}=(a+b\sqrt{2})(c+d\sqrt{2})$ for $r\in\{1,2\}$ and $a,b,c,d\in\nn_0$, then $ac+2bd=1$. Hence $a=c=1$ and $bd=0$, so one of the factors is $1$. Hence $1+\sqrt{2}$ and $1+2\sqrt{2}$ are atoms. Suppose now that $19+3\sqrt{2}$ factors in the same way. Comparing coefficients, we find that $ac+2bd=19$ and $ad+bc=3$. In particular, $a,c>0$. If either $b$ or $d$ is zero then the corresponding integer factor divides both $19$ and $3$ and must be $1$. Otherwise $ad$ and $bc$ are positive integers with sum $3$, so $abcd=2$. But then $ac+2bd$ is either $4$ or $5$, a contradiction. Thus $19+3\sqrt{2}$ is also an atom. The equalities
    \[
        (1+2\sqrt{2})^3=25+22\sqrt{2}=(19+3\sqrt{2})(1+\sqrt{2})
    \]
    exhibit two factorizations of the same element of $S^\bullet$ into atoms of different lengths.
\end{example}

We turn to the transfer Krull property. Recall that every half-factorial monoid is transfer Krull, but such a monoid need not be root-closed. For instance, let $M$ be the submonoid of $(\nn_0^2,+)$ generated by $(1,0)$, $(1,2)$, and $(1,3)$. These three generators are the atoms of $M$, and the length of any factorization of an element of $M$ is its first coordinate. Therefore $M$ is half-factorial. However, $(1,1)\notin M$, although $(1,1)=(1,3)-(1,2)+(1,0)$ belongs to $\gp(M)$ and $2(1,1)=(1,0)+(1,2)\in M$. Hence transfer Krullness alone does not imply the root closedness used in the proof of Theorem~\ref{thm:every-element-algebraic}.

For each transcendental real number $t>0$, the additive monoid of $\nn_0[t]$ is free on $\{t^n:n\in\nn_0\}$. It would be natural to look for a positive semiring of this form whose multiplicative monoid is transfer Krull. To see why this is not possible, we use the polynomial of Kainrath from Section~\ref{sec:background}. Recall that in a root-closed monoid no atom can be a proper power of an element of its Grothendieck group. The following lemma extends this observation to transfer Krull monoids.

\begin{lemma}\label{lem:power-obstruction}
    Let $M$ be a transfer Krull monoid, and let $u,v,w\in M$ and $k\in\nn$ with $k\ge 2$ satisfy $u^k=vw^k$. Then $v$ is not an atom of $M$.
\end{lemma}

\begin{proof}
    Take a transfer homomorphism $\theta\colon M\to M_0$ into a reduced Krull monoid and a divisor theory $\partial\colon M_0\to\nn_0^{(P)}$. Applying $\partial\circ\theta$ to the equality $u^k=vw^k$ yields
    \[
        k\,\partial(\theta(u))=\partial(\theta(v))+k\,\partial(\theta(w)).
    \]
    In particular, $\partial(\theta(w))\le\partial(\theta(u))$, whence $\theta(u)=\theta(w)y$ for some $y\in M_0$. Substituting this expression for $\theta(u)$ into the displayed equality, we see that $\partial(\theta(v))=k\,\partial(y)=\partial(y^k)$. The injectivity of $\partial$ implies that $\theta(v)=y^k$. Since $k\ge 2$, the element $\theta(v)$ is not an atom of $M_0$, and therefore $v$ is not an atom of $M$.
\end{proof}

\begin{example}\label{ex:polynomial-semiring}
    It was established in~\cite[Remark~5.4]{CF19} that $\nn_0[x]^\bullet$ is not transfer Krull. We present another proof using Lemma~\ref{lem:power-obstruction}. Consider the polynomial $f=x^4+2x^3-x^2+4x+2$ of Kainrath. By~\eqref{eq:kainrath-multiple} and~\eqref{eq:kainrath-square}, both $(x+1)f$ and $f^2$ belong to $\nn_0[x]^\bullet$. By Lemma~\ref{lem:power-obstruction} and the equality $((x+1)f)^2=f^2(x+1)^2$, it suffices to show that $f^2$ is an atom of $\nn_0[x]^\bullet$.

    First observe that $f$ is irreducible in $\zz[x]$. Its only possible rational roots are $\pm1$ and $\pm2$, and none of these is a root. If $f=(x^2+ax+b)(x^2+cx+d)$ for some $a,b,c,d\in\zz$, then comparing coefficients yields $bd=2$, $a+c=2$, and $ac=-1-b-d$. Since $bd=2$, we have $b+d\in\{-3,3\}$, so $ac\in\{-4,2\}$. Neither value is possible for integers $a,c$ with sum $2$. Thus $f$ is irreducible in $\zz[x]$. Because neither $f$ nor $-f$ has nonnegative coefficients, the factoriality of $\zz[x]$ implies that $f^2$ is an atom of $\nn_0[x]^\bullet$. Consequently, $\nn_0[x]^\bullet$ is not transfer Krull.
\end{example}

In particular, if $t>0$ is transcendental then $\nn_0[t]$ cannot have a transfer Krull multiplicative monoid, even though its additive monoid is factorial. We conclude with the following question.

\begin{question}\label{ques:mixed}
    Let $S$ be a positive semiring such that $(S,+)$ is Krull and $S^\bullet$ is transfer Krull. Must $S=\nn_0$?
\end{question}

\bigskip
\section*{Acknowledgements}

The authors participated in the CrowdMath Internship at MIT during the preparation of this paper. They thank the program for making this collaboration possible. The fourth author was supported by a National Science Foundation Postdoctoral Fellowship.

\bigskip

\end{document}